\documentclass[12pt]{article}
\usepackage{amsmath,amssymb,amsfonts}

\newtheorem{theorem}{Theorem}
\newtheorem{proposition}{Proposition}
\newtheorem{corollary}{Corollary}

\def\D{{\cal D}}
\def\A{{\cal A}}
\def\B{{\cal B}}
\def\R{{\mathbb R}}
\def\C{{\mathbb C}}

\def\n{{\bf n}}
\def\x{{\bf x}}
\def\X{{\bf X}}

\begin{document}

\title{The Moutard transformation of two-dimensional Dirac operators and the~M\"obius geometry}
\author{Iskander A. TAIMANOV
\thanks{Sobolev Institute of Mathematics, Academician Koptyug avenue 4, 630090, Novosibirsk, Russia, and Department of Mathematics and Mechanics, Novosibirsk State University, Pirogov street 2, 630090 Novosibirsk, Russia; e-mail: taimanov@math.nsc.ru. \newline
The work was supported by RSF (grant 14-11-00441). }}
\date{}
\maketitle

\section{Introduction}

In the present article we consider the Moutard transformation for two-di\-men\-sio\-nal Dirac operators
\begin{equation}
\label{dirac}
\D = \left(
\begin{array}{cc}
0 & \partial \\
-\bar{\partial} & 0
\end{array}
\right) + \left(
\begin{array}{cc}
U & 0 \\
0 & U
\end{array}
\right),
\end{equation}
where $\partial = \frac{1}{2}\big(\frac{\partial}{\partial x} - i\frac{\partial}{\partial y}\big)$ and
$\bar{\partial} =
\frac{1}{2}\big(\frac{\partial}{\partial x} + i\frac{\partial}{\partial y}\big)$,
with real-valued potentials $U$. This transformation and its extension to a transformation of solutions to
the modified Novikov--Veselov equation was introduced in \cite{C}.

The Weierstrass representation of surfaces corresponds to a solution
$$
\psi =
\left(\begin{array}{c}
\psi_1 \\ \psi_2
\end{array}
\right)
$$
of the Dirac equation
\begin{equation}
\label{diraceq0}
\D \psi = 0,
\end{equation}
a surface in $\R^3$ determined by the classical Weierstrass formulas for
minimal surfaces in $\R^3$ \cite{K} and moreover any surface in $\R^3$ admits such
a representation even globally (wherewith $\psi$ is a section of spinor bundle over a surface)
and important geometrical properties of the surface (in particular, the value of Willmore functional) are encoded into
the spectral characteristics of $\D$  \cite{T-MNV,T-RS}.

Hence to every surface $S \subset \R^3$ with a fixed conformal parameter $z=x+iy$
there corresponds a unique Dirac operator $\D$, coming into its Weierstrass representation, with the potential
$U$ (the potential of a surface with a fixed conformal parameter).

In the present article we

\begin{enumerate}
\item
{\sl describe the action of the inversion on the data of the Weierstrass representation, i.e. on $U$ and $\psi$
(Theorem 1)};

\item
{\sl show that the Moutard transformation} \cite{C} {\sl has a geometrical meaning:
the Moutard transformation maps the potential $U$ of a surface  $S$ into the potential $\widetilde{U}$ of
the inverted surface $\widetilde{S}$
(Theorem 2).}
\end{enumerate}

\section{Preliminary facts}
\label{section2}

\subsection{The Weierstrass representation of surfaces}

Let
$$
r: {\cal U} \to \R^3, \ \ \ {\cal U} \subset \C,
$$
be a conformal immersion of a domain ${\cal U}$ into $\R^3$, i.e., that is an immersion such that the induced metric takes the form
\begin{equation}
\label{induced}
ds^2 = e^{2\alpha(z,\bar{z})} dz\, d\bar{z}.
\end{equation}
Such a parameter $z = x+iy$ is called a conformal parameter on the surface.
The conformality condition reads
$$
\left(\frac{\partial x^1}{\partial z}\right)^2 +
\left(\frac{\partial x^2}{\partial z}\right)^2 +
\left(\frac{\partial x^3}{\partial z}\right)^2  = 0,
$$
where $x^1,x^2,x^3$ are the Euclidean coordinates in $\R^3$. The points of this quadric in $\C^3$
are parameterized by pairs $(\psi_1,\bar{\psi}_2) \in \C^2$ as follows:
$$
\frac{\partial x^1}{\partial z} = \frac{i}{2} (\psi_1^2 + \bar{\psi}_2^2), \ \ \
\frac{\partial x^2}{\partial z} = \frac{1}{2} (\bar{\psi}_2^2 - \psi_1^2), \ \ \
\frac{\partial x^3}{\partial z} = \psi_1 \bar{\psi}_2.
$$
The function $\psi = \left(\begin{array}{c} \psi_1 \\ \psi_2 \end{array}\right)$
defines a surface via these formulas if and only if
it satisfies the Dirac equation
(\ref{diraceq0}) for the Dirac operator (\ref{dirac}) with a real-valued potential $U$.
The Dirac equation reads exactly that the form
$$
\sum_{k=1}^3 \left(\frac{\partial x^k}{\partial z} dz + \overline{\frac{\partial x^k}{\partial z}} d\bar{z}\right)
$$
is closed and the coordinates $x^1,x^2,x^3$ are real-valued.

Therewith, given any solution $\psi$ of (\ref{diraceq0}),
we construct a surface $r: {\cal U} \to \R^3$ by formulas
$$
x^1(P) = \frac{i}{2} \int_{P_0}^P
\left((\psi_1^2 + \bar{\psi}_2^2)dz - (\bar{\psi}_1^2 + \psi_2^2)d\bar{z}\right) + x^1(P_0),
$$
\begin{equation}
\label{weier}
x^2(P) = \frac{1}{2} \int_{P_0}^P
\left((\bar{\psi}_2^2 - \psi_1^2)dz + (\psi_2^2 - \bar{\psi}_1^2)d\bar{z}\right) + x^2(P_0),
\end{equation}
$$
x^3(P) = \int_{P_0}^P
\left(\psi_1 \bar{\psi}_2 dz + \bar{\psi}_1 \psi_2 \bar{z}\right) + x^3(P_0).
$$
Here $P_0$ is a fixed point in ${\cal U}$ and the integral is taken over a path in ${\cal U}$ joining $P_0$ and $P$.
If ${\cal U}$ is simply-connected, then the integral does not depend on a path.
The induced metric (\ref{induced}) takes the form
$$
e^\alpha = |\psi_1|^2 + |\psi_2|^2,
$$
and hence the formulas (\ref{weier}) define an immersion of ${\cal U}$ exactly outside branch points at which $|\psi_1|^2 = |\psi_2|^2=0$.
The unit normal vector $\n$ equals to
\begin{equation}
\label{normal}
\n = e^{-\alpha} ( i(\psi_1\psi_2 -\bar{\psi}_1\bar{\psi}_2), - (\psi_1\psi_2 +\bar{\psi}_1\bar{\psi}_2),
(|\psi_2|^2 - |\psi_1|^2)),
\end{equation}
and the potential $U$ of the Dirac operator is equal to
\begin{equation}
\label{potential}
U = \frac{e^\alpha H}{2},
\end{equation}
where $H$ is the mean curvature of the surface.

Since (\ref{weier}) define a surface up to translations, the data $U$ and $\psi$ are invariant under them. Rotations of
$\R^3$ preserve $U$ and induce spinor actions on $\psi$  \cite{T-WN}.

\subsection{The inversion}

The M\"obius geometry studies geometrical figures in the Euclidean spaces complemented by the infinity points, i.e. in the spheres $S^n = \R^n \cup \{\infty\}$,
and their properties which are invariant with respect to conformal transformations.

For $n \geq 3$ the group formed by all orientation preserving conformal transformations is generated by translations and rotations of $\R^n$ and the inversion
and is isomorphic to the unity component of
$SO(n+1,1)$.

The inversion of the three-space is as follows
$$
T: \x \to -\frac{\x}{|\x|^2}, \ \ \ \x = (x^1,x^2,x^3) \in \R^3.
$$
It maps conformally $\R^3 \cup \{\infty\}$ onto itself.

Let $u$ be a vector tangent to $\R^3$ at $\x$: $u \in T_{\x}\R^3$ and let $\x\neq 0$. Then it is easy to compute that
\begin{equation}
\label{inversionvector}
T^\ast u =  -\frac{u}{|\x|^2} + 2\x \frac{\langle \x,u\rangle}{|\x|^4},
\end{equation}
where $\langle \cdot,\cdot \rangle$ is the Euclidean scalar product.
This implies that
\begin{equation}
\label{inversionproduct}
\langle T^\ast u, T^\ast v\rangle = \frac{\langle u,v \rangle}{|\x|^4}, \ \ \ \ u,v \in T_{\x}\R^3.
\end{equation}

Let
$$
r: {\cal U} \to \R^3
$$
be a surface with a conformal parameter $z$.
The inversion transforms it into the surface:
$$
\widetilde{r} = T \cdot r: {\cal U} \to \R^3.
$$
By (\ref{inversionproduct}), $z$ is a conformal parameter for $\widetilde{r}$ and the conformal factors are
related by the formula:
\begin{equation}
\label{factor}
e^{\widetilde{\alpha}(z,\bar{z})} = \frac{e^{\alpha(z,\bar{z})}}{|r(z,\bar{z})|^2}, \ \ \ e^{\widetilde{\alpha}} =
\frac{1}{2} \langle \widetilde{r}_z, \widetilde{r}_{\bar{z}}\rangle, \ \ e^\alpha = \frac{1}{2} \langle r_z, r_{\bar{z}}\rangle.
\end{equation}

Let us recall that
$$
\Delta r = 4 e^{-2\alpha} \partial \bar{\partial} r = 2H\n,
$$
which we rewrite as
\begin{equation}
\label{delta}
r_{z\bar{z}} = e^\alpha U \n,
\end{equation}
where $\n$ is the unit normal vector to the surface and $U$ is the potential of the Weierstrass representation of $r$.
The latter formula implies that
$$
U = e^{-\alpha} \langle r_{z\bar{z}}, \n \rangle.
$$

\begin{proposition}
The potential $\widetilde{U}$ of the Weierstrass representation of the surface $\widetilde{r}$ (with a conformal parameter $z$)
is equal to
\begin{equation}
\label{tildepotential}
\widetilde{U} = U + e^\alpha \frac{\langle r, \n \rangle}{|r|^2}.
\end{equation}
\end{proposition}

{\sc Proof.}
By (\ref{inversionvector}) and (\ref{inversionproduct}), we have
$$
\widetilde{r}_z = -\frac{r_z}{|r|^2} + \frac{2r\langle r,r_z\rangle}{|r|^4}, \ \ \ \
\widetilde{\n} = -\n + \frac{2r\langle r,\n\rangle}{|r|^2}.
$$
By straightforward commutations we derive
$$
\widetilde{r}_{z\bar{z}} = -\frac{r_{z\bar{z}}}{|r|^2}  - \frac{8r\langle r,r_z\rangle \langle r, r_{\bar{z}}\rangle}{|r|^6} +
$$
$$
+
 \frac{2}{|r|^4} \left( r_z \langle r, r_{\bar{z}} \rangle + r_{\bar{z}} \langle r, r_z \rangle + r \langle r_z, r_{\bar{z}} \rangle
+ r \langle r, r_{z\bar{z}} \rangle \right)
$$
and
$$
\langle \widetilde{r}_{z\bar{z}}, \widetilde{\n} \rangle = \frac{\langle r_{z\bar{z}}, \n \rangle}{|r|^2} +
2\frac{\langle r, \n \rangle \langle r_z, r_{\bar{z}} \rangle}{|r|^4}.
$$
By (\ref{delta}), we rewrite the latter equation as follows
$$
e^{\widetilde{\alpha}} \widetilde{U} = \frac{e^\alpha U}{|r|^2} + \frac{\langle r, \n\rangle e^{2\alpha}}{|r|^4}
$$
and after dividing both sides by $e^{\widetilde{\alpha}} = e^\alpha/|r|^2$ (see (\ref{factor})) we obtain (\ref{tildepotential}).
Proposition is proved.

\subsection{The modified Novikov--Veselov equation}

The modified Novikov--Veselov equation (the mNV equation)
$$
U_t = \big(U_{zzz} + 3U_z V + \frac{3}{2}UV_z \big) + \big(U_{\bar{z}\bar{z}\bar{z}} + 3U_{\bar{z}}\bar{V} + \frac{3}{2} U\bar{V}_{\bar{z}}\big),
$$
where
$$
V_{\bar{z}} = (U^2)_z,
$$
was introduced
in \cite{Bogdanov}.

This equation takes the form of Manakov triple:
$$
\D_t + [\D,\A] - \B\D = 0,
$$
where $\D$ is a two-dimensional Dirac operator (\ref{dirac}).
Usually the Manakov representation of the mNV equation
$$
L_t  + [L,A]-BL = 0
$$
was written in terms of the operator $L$ of the form
$$
L = \left(
\begin{array}{cc}
 \partial & -U \\
U & \bar{\partial}
\end{array}
\right)
$$
(see \cite{C,T-MNV}).
Since
$$
\D = L \Gamma
$$
with
\begin{equation}
\label{gamma}
\Gamma = \left(\begin{array}{cc} 0 & 1 \\ -1 & 0 \end{array}\right),
\end{equation}
these representation are related by the formulas
$$
\A = - \Gamma A \Gamma, \ \ \ \B = \Gamma A \Gamma + A + B,
$$
with
\begin{equation}
\label{a}
\A = \partial^3 + \bar{\partial}^3 +
\end{equation}
$$
+ 3\left(\begin{array}{cc}  V & 0 \\ U_z & 0 \end{array}\right)\partial +
3\left(\begin{array}{cc} 0 & -U_{\bar{z} }\\ 0 & \bar{V} \end{array}\right) \bar{\partial}+
\frac{3}{2}\left(\begin{array}{cc} V_z & 2U\bar{V} \\ -2UV & \bar{V}_{\bar{z}} \end{array}\right),
$$
$$
\B =
3\left(\begin{array}{cc} -V & 0 \\ -2U_z & V \end{array}\right)\partial +
3\left(\begin{array}{cc} \bar{V} & 2U_{\bar{z}} \\ 0 & -\bar{V} \end{array}\right)\bar{\partial} +
\frac{3}{2}\left( \begin{array}{cc} \bar{V}_{\bar{z}} - V_z & 2 U_{\bar{z}\bar{z}} \\ -2U_{zz} & V_z - \bar{V}_{\bar{z}}\end{array} \right)
$$

\subsection{The matrix algebra $H$}

We note that
$$
\Gamma = i\sigma_2,
$$
where $\sigma_2$ is one of the Pauli matrices which are
$$
\sigma_1 = \left(
\begin{array}{cc}
0 & 1 \\
1 & 0
\end{array}
\right), \ \ \
\sigma_2 = \left(
\begin{array}{cc}
0 & -i \\
i & 0
\end{array}
\right), \ \ \
\sigma_3 = \left(
\begin{array}{cc}
1 & 0 \\
0 & -1
\end{array}
\right).
$$
The Pauli matrices satisfy the relations
$$
\sigma_a \sigma_b = i\varepsilon_{abc} \sigma_c + \delta_{ab}\sigma_0,
$$
where $\sigma_0$ is the unity matrix, and  $i\sigma_1, i\sigma_2,i\sigma_3$ form a basis for the Lie algebra
$su(2)$ formed by all matrices of the form  $A = \left(\begin{array}{cc} a & b \\ -\bar{b} & \bar{a} \end{array}\right), \mathrm{Tr}\,A=0$.

The four-dimensional space $H$ formed by all matrices of the form
$$
\left(\begin{array}{cc} a & b \\ -\bar{b} & \bar{a} \end{array}\right), \ \ a,b \in \C,
$$
is spanned over $\R$ by $\sigma_0$, $i\sigma_1, i\sigma_2,i\sigma_3$,
is isomorphic to the quaternion algebra. In particular, $H$ is closed with respect to the product.

\section{The Moutard transformation}
\label{section3}

In this section we expose the Moutard type transformation for two-dimensional Dirac operators and solutions of the mNV equation introduced in \cite{C}.
However we modify the initial presentation for demonstrating the geometry which is  hidden
in ana\-ly\-tical formulas and was unnoticed until recently.

\subsection{The Moutard transformation of operators}

If a spinor
$\psi = \left(\begin{array}{c} \psi_1 \\ \psi_2 \end{array}\right)$ which meets the Dirac equation
(\ref{diraceq0}), then it is easy to notice that
$$
\psi^\ast = \left(\begin{array}{c} -\bar{\psi}_2 \\ \bar{\psi}_1 \end{array}\right)
$$
also satisfies this equation. Let us form  from $\psi$ and $\psi^\ast$ a matrix-valued solution
\begin{equation}
\label{matrixeq}
\Psi =
\left(\begin{array}{cc}  \psi_1 & -\bar{\psi}_2 \\  \psi_2 & \bar{\psi}_1 \end{array}\right)
\end{equation}
of the Dirac equation (\ref{diraceq0}).
We note that $\Psi$ is an $H$-valued function.

For every pair $\Psi$ and
$\Phi$ of $H$-valued functions
let us  correspond a matrix-valued $1$-form $\omega$
\begin{equation}
\label{omega}
\omega(\Phi,\Psi) = \Phi^\top \Psi dy - i \Phi^\top \sigma_3 \Psi dx
=
\end{equation}
$$
-\frac{i}{2}\left(\Phi^\top \sigma_3 \Psi + \Phi^\top \Psi\right) dz - \frac{i}{2}\left(\Phi^\top \sigma_3 \Psi - \Phi^\top \Psi\right) d\bar{z}
$$
and a matrix-valued function
\begin{equation}
\label{sigma}
S(\Phi,\Psi)(z,\bar{z},t) = \Gamma \int_0^z \omega(\Phi,\Psi),
\end{equation}
which  is defined up to constant matrices from $su(2)$ formed by integration constants.

Here and in the sequel we denote the transposition of $X$ by $X^\top$.

To every $H$-valued function $\Psi$ we correspond a matrix-valued function
\begin{equation}
\label{kmatrix}
K(\Psi) =  \Psi S^{-1}(\Psi,\Psi)\Gamma \Psi^\top\Gamma^{-1}.
\end{equation}

\begin{proposition}
[\cite{C}]
\label{prop2}
Given a solution $\Psi_0$ of the Dirac equation (\ref{matrixeq}) for the operator $\D$ with real-valued potential $U$,
the matrix $K(\Psi_0)$ takes the form
$$
K(\Psi_0) =
\left(\begin{array}{cc} iW & a \\ -\bar{a} & -iW
\end{array}\right)
$$
with $W$ real-valued, and
for every solution $\Psi$ of form (\ref{matrixeq}) of the Dirac equation
(\ref{diraceq0}) the function $\widetilde{\Psi}$ of the form
\begin{equation}
\label{moutard1}
\widetilde{\Psi} =  \Psi - \Psi_0 S^{-1}(\Psi_0,\Psi_0) S(\Psi_0,\Psi)
\end{equation}
satisfies the equation
$$
\widetilde{\D}\widetilde{\Psi} = 0
$$
for the Dirac operator $\widetilde{\D}$ with potential
$$
\widetilde{U} = U + W.
$$
\end{proposition}

\subsection{The Moutard transformation of solutions to the mNV equation}

Given a pair of $H$-valued functions $\Phi$ and $\Psi$ which depend on $x,y$, and $t$, let us define a  matrix-valued $1$-form
$$
\widehat{\omega}(\Phi,\Psi) = \Phi^\top \Psi dy - i \Phi^\top \sigma_3 \Psi dx +
\left[i(\Phi^\top_{yy}\sigma_3\Psi + \Phi^\top\sigma_3\Psi_{yy}-\Phi^\top_y\sigma_3\Psi_y)
\right. +
$$
\begin{equation}
\label{omega1}
\left.
2iU(\Phi^\top_y\sigma_2\Psi - \Phi^\top\sigma_2\Psi_y)
+
\Phi^\top
\left(\begin{array}{cc}
iU^2-3iV & -iU_x \\
-iU_x & -iU^2 +3i\bar{V}
\end{array}\right)
\Psi
\right]dt =
\end{equation}
$$
\omega(\Phi,\Psi) +
\left[-i((\Phi^\top_{zz} + \Phi^\top_{\bar{z}\bar{z}} - 2\Phi^\top_{z\bar{z}})\sigma_3\Psi +
\Phi^\top\sigma_3(\Psi_{zz} + \Psi_{\bar{z}\bar{z}} - 2 \Psi_{z\bar{z}}) -
\right.
$$
$$
(\Phi^\top_z - \Phi^\top_{\bar{z}})\sigma_3 (\Psi_z - \Psi_{\bar{z}})) -
2U((\Phi^\top_z - \Phi^\top_{\bar{z}}) \sigma_2 \Psi - \Phi^\top\sigma_2(\Psi_z - \Psi_{\bar{z}}))
+
$$
$$
\left.
\Phi^\top
\left(\begin{array}{cc}
iU^2-3iV & -i(U_z + U_{\bar{z}}) \\
-i(U_z + U_{\bar{z}}) & -iU^2 +3i\bar{V}
\end{array}\right)
\Psi
\right]dt,
$$
and matrix-valued functions
$$
\widetilde{S}(\Phi,\Psi)(z,\bar{z},t) = \Gamma \int_0^z \widetilde{\omega}(\Phi,\Psi),
$$
$$
K(\Psi) =  \Psi \widetilde{S}^{-1}(\Psi,\Psi)\Gamma \Psi^\top\Gamma^{-1}, \ \ \
$$
$$
M(\Psi) = \Gamma \Psi_y \Psi^{-1} \Gamma^{-1} =i\Gamma (\Psi_z - \Psi_{\bar{z}})\Psi^{-1}\Gamma^{-1}.
$$

\begin{proposition}
[\cite{C}]
\label{prop3}
Let $U(z,\bar{z},t)$ and $V(z,\bar{z},t)$ satisfy the mNV equation and $\Psi_0(z,\bar{z},t)$ satisfy the system
$$
\D \Psi_0 = 0, \ \ \ \ \frac{\partial \Psi_0}{\partial t } = \A \Psi_0.
$$
Then

\begin{enumerate}
\item
the matrices $K(\Psi_0)$ and $M(\Psi_0)$ take the form
  $$
K =
\left(\begin{array}{cc} iW & a \\ -\bar{a} & -iW
\end{array}\right), \ \ \
M =
\left(\begin{array}{cc} b & c \\ -\bar{c} & \bar{b}
\end{array}\right),
$$
with $W$ real valued;

\item
for every solution $\Psi$ (\ref{matrixeq}) of the Dirac equation (\ref{diraceq0}) and
$$
\frac{\partial \Psi}{\partial t}= \A \Psi
$$
the function $\widetilde{\Psi}$ of the form
$$
\widetilde{\Psi} =  \Psi - \Psi_0 \widetilde{S}^{-1}(\Psi_0,\Psi_0) \widetilde{S}(\Psi_0,\Psi)
$$
satisfies the equation
$$
\widetilde{\D}\widetilde{\Psi} = 0
$$
for the Dirac operator $\widetilde{\D}$ with potential
$$
\widetilde{U} =  U + W
$$
and the equation
$$
\frac{\partial \widetilde{\Psi}}{\partial t} = \widetilde{\A} \widetilde{\Psi}
$$
where $\widetilde{\A}$ takes the form (\ref{a}) with $U$ replaced by $\widetilde{U}$ and $V$ replaced by
$\widetilde{V}$:
$$
\widetilde{V} =  V + 2UW  + a^2 + 2(a\bar{b} - i\bar{c}W);
$$

\item
the function $\widetilde{U}$ is real-valued and $\widetilde{U}$ and $\widetilde{V}$ satisfy the mNV equation
$$
\widetilde{U}_t = \big(\widetilde{U}_{zzz} + 3\widetilde{U}_z \widetilde{V} +
\frac{3}{2}\widetilde{U}\widetilde{V}_z \big) + \big(\widetilde{U}_{\bar{z}\bar{z}\bar{z}} +
3\widetilde{U}_{\bar{z}}\bar{\widetilde{V}} + \frac{3}{2} \widetilde{U}\bar{\widetilde{V}}_{\bar{z}}\big),
$$
$$
\widetilde{V}_{\bar{z}} = (\widetilde{U}^2)_z.
$$
\end{enumerate}
\end{proposition}

\section{The action of the inversion on the Weierstrass representation data}
\label{section4}

Let $\psi = \left(\begin{array}{c} \psi_1 \\ \psi_2 \end{array}\right)$ define a surface
$r:  {\cal U} \to \R^3$
via the Weierstrass representation.

Let us identify $\R^3$ with the Lie algebra $su(2)$ via the mapping
$$
\x = (x^1,x^2,x^3) \to
\X = \left(\begin{array}{cc} ix^3 & -x^1 - ix^2 \\ x^1-ix^2 & -ix^3 \end{array}\right) .
$$
In such a representation the inversion takes a very simple form:
$$
T(\X) = \X^{-1}.
$$

Let us construct a matrix-valued function
$$
\Psi_0 = \left(\begin{array}{cc} \psi_1 & -\bar{\psi}_2 \\ \psi_2 & \bar{\psi}_1 \end{array}\right).
$$

\begin{proposition}
The formula (\ref{sigma}) gives an immersion into $su(2) = \R^3$ of a surface defined by the spinor
$\psi$ via the Weierstras representation.
\end{proposition}

{\sc Proof.}
By (\ref{weier}) and (\ref{omega}),
we have
$$
S(\Psi_0,\Psi_0)(P) = \Gamma \int_{P_0}^P
-\frac{i}{2}\left(\Psi_0^\top (\sigma_3+1) \Psi_0 dz + \Psi_0^\top (\sigma_3 -1)
\Psi_0\right) d\bar{z}) =
$$
\begin{equation}
\label{surface}
= i \int_0^P \left(\begin{array}{cc} \psi_1 \bar{\psi}_2 & -\bar{\psi}_2^2 \\
\psi_1^2 & -\psi_1\bar{\psi}_2 \end{array}\right)dz +
\left(\begin{array}{cc} \bar{\psi}_1 \psi_2 & \bar{\psi}_1^2 \\
-\psi_2^2 & -\bar{\psi}_1 \psi_2 \end{array}\right)d\bar{z} =
\end{equation}
$$
= \int_0^P d \left(\begin{array}{cc} ix^3 & -x^1 - ix^2 \\ x^1-ix^2 & -ix^3 \end{array}\right) \in su(2) ,
$$
i.e. $S$ determines  a surface defined by $\psi$ via the Weierstrass representation.
Proposition is proved.

Let us fix integration constants to achieve
$$
S(\Psi_0,\Psi_0)(P) = r(P).
$$

Let us consider the inversion $
\widetilde{r} = T \cdot r: {\cal U} \to \R^3,
$ of the surface $r$.
It is given by the formula
$$
P \in {\cal U} \to \widetilde{S}(P) = S^{-1}(P).
$$
We have
$\widetilde{S}(P) = S(\widetilde{\Psi}_0,\widetilde{\Psi}_0)(P)  = S^{-1}(\Psi_0,\Psi_0)(P)$,
and, by (\ref{omega}) and (\ref{sigma}),
$$
\widetilde{S} _z= -\frac{i}{2} \Gamma \widetilde{\Psi}_0^\top (1+\sigma_3)\widetilde{\Psi}_0.
$$
Since $\widetilde{S}S = S^{-1}S = 1$,
we have
$(\widetilde{S}S)_z = \widetilde{S}_z S + \widetilde{S} S_z = 0$
and therefore conclude that
$$
\widetilde{S}_z = - S^{-1} S_z S^{-1}.
$$
This implies the equality
$$
-\frac{i}{2} \Gamma \widetilde{\Psi}_0^\top (1+\sigma_3)\widetilde{\Psi}_0 =
\frac{i}{2} S^{-1} \Gamma \Psi_0^\top (1+ \sigma_3) \Psi_0 S^{-1}
$$
which after simple cancellations takes the form
$$
\widetilde{\Psi}_0^\top (1+\sigma_3)\widetilde{\Psi}_0 =
-\Gamma^{-1} S^{-1} \Gamma \Psi_0^\top (1+ \sigma_3) \Psi_0 S^{-1}.
$$
Since $S^{-1}$ equals to
$$
S^{-1} = \frac{1}{\det S} \left(\begin{array}{cc} - ix^3 & x^1 + i x^2 \\ -x^1 + ix^2 & ix^3 \end{array}\right),
$$
we have
$$
-\Gamma^{-1} S^{-1}\Gamma = \Gamma S^{-1} \Gamma = (S^{-1})^\top
$$
which implies
\begin{equation}
\label{z}
C^\top (1+\sigma_3) C = (1+\sigma_3) = \left(\begin{array}{cc} 2 & 0 \\ 0 & 0 \end{array}\right)
\end{equation}
with
$$
C = \Psi_0 S^{-1} \widetilde{\Psi}_0^{-1}.
$$
Analogously by considering $\widetilde{S}_{\bar{z}}$ and $S_{\bar{z}}$ we derive that
$$
C^\top (\sigma_3-1)C = (\sigma_3-1) =  \left(\begin{array}{cc} 0 & 0 \\ 0 & -2 \end{array}\right).
$$
It follows from the latter equality together with (\ref{z}) that $C$ is diagonal and $C^\top = C^{-1}$.
Since $\Psi_0, S^{-1}, \widetilde{\Psi}_0 \in H$,  we conclude that $C\in H$ which implies that
$C = \pm \left(\begin{array}{cc} 1 & 0 \\ 0 & 1 \end{array}\right)$, and hence
$\widetilde{\Psi}_0 = \pm \Psi_0 S^{-1}$.
Since both spinors $\pm \Psi_0 S^{-1}$ define the same surface, we put without loss of generality
$$
\widetilde{\Psi}_0 = \Psi_0 S^{-1}.
$$

The spinor $\widetilde{\Psi}_0$ satisfies the Dirac equation
$$
\widetilde{D} \widetilde{\Psi}_0 = (\D_0 + \widetilde{U})\widetilde{\Psi}_0 = 0
$$
where
$$
\D_0 = \left(\begin{array}{cc} 0 & \partial \\ -\bar{\partial} & 0 \end{array}\right).
$$
It is easy to check the ``Leibniz rule''
\begin{equation}
\label{leibniz}
\D_0 (A \cdot B) = (\D_0\,A) \cdot B +
\left(\begin{array} {cc} 0 & 1 \\ 0 & 0 \end{array}\right) A \cdot  \partial B
+
\left(\begin{array} {cc} 0 & 0 \\ -1 & 0 \end{array}\right) A \cdot \bar{\partial} B
\end{equation}
and apply it as follows
$$
\D_0 (\Psi_0 S^{-1}) = (\D_0 \Psi_0)S^{-1} +
\left(\begin{array} {cc} 0 & 1 \\ 0 & 0 \end{array}\right) \Psi_0 S^{-1}_z
+ \left(\begin{array} {cc} 0 & 0 \\ -1 & 0 \end{array}\right) \Psi_0 S^{-1}_{\bar z} =
$$
$$
= - U\Psi_0 S^{-1}  +
 i \left(\begin{array} {cc} 0 & 1 \\ 0 & 0 \end{array}\right) \Psi_0 S^{-1} \Gamma \Psi_0^\top
 \left(\begin{array}{cc} 1 & 0 \\ 0 & 0 \end{array}\right)\Psi_0  S^{-1} +
 $$
 $$
  + i \left(\begin{array} {cc} 0 & 0 \\ -1 & 0 \end{array}\right)  \Psi_0 S^{-1}\Gamma \Psi_0^\top
\left(\begin{array} {cc} 0 & 0 \\ 0 & -1 \end{array}\right) \Psi_0 S^{-1}.
$$
By dividing both sides by $\Psi_0 S^{-1}$, we derive
\begin{equation}
\label{inversionu}
\widetilde{U} = U -
i \left(\left(\begin{array} {cc} 0 & 1 \\ 0 & 0 \end{array}\right) G
 \left(\begin{array}{cc} 1 & 0 \\ 0 & 0 \end{array}\right)
+ \left(\begin{array} {cc} 0 & 0 \\ 1 & 0 \end{array}\right) G
\left(\begin{array} {cc} 0 & 0 \\ 0 & 1 \end{array}\right)\right),
\end{equation}
$$
G =  \Psi_0 S^{-1} \Gamma \Psi_0^\top.
$$
Finally we conclude

\begin{theorem}
Let $\Psi_0$ define a surface $S$ of the form (\ref{surface}) via the Weierstrass representation and
meet the Dirac equation with potential $U$.
Let $\widetilde{S}$ be a surface obtained from $S$ by the inversion.

Then
$\widetilde{S}$ is defined by the spinor
$$
\widetilde{\Psi}_0 = \Psi_0 S^{-1}
$$
via the Weierstrass
representation and $\widetilde{\Psi}_0$ meets the Dirac equation with the potential $\widetilde{U}$
of the form (\ref{inversionu}).
\end{theorem}

\section{Geometry of the Moutard transformation}
\label{section5}

\begin{theorem}
\label{th2}
The Moutard transformation of the Dirac operator $\D$ given in Proposition \ref{prop2} maps the potential
of the Weierstrass representation of $S$ (\ref{surface}) into the potential of its inversion $\widetilde{S}$.
\end{theorem}

{\sc Proof.}
Let us compute $W=-iK_{11}$ given by (\ref{kmatrix}). We have
$$
K = \Psi_0 S^{-1} \Gamma \Psi_0^\top \Gamma^{-1} =
\Psi_0 S^{-1} \left(\begin{array}{cc} 0 & 1 \\ -1 & 0 \end{array}\right)
\left(\begin{array}{cc} \psi_1 & \psi_2 \\ -\bar{\psi}_2 & \bar{\psi}_1 \end{array}\right)
\left(\begin{array}{cc} 0 & -1 \\ 1 & 0 \end{array}\right) =
$$
$$
=
\frac{1}{|r|^2}
\left(\begin{array}{cc} \psi_1 & -\bar{\psi}_2 \\ \psi_2 & \bar{\psi}_1 \end{array}\right)
\left(\begin{array}{cc} -ix^3 & x^1 + ix^2 \\ -x^1 +ix^2 & ix^3 \end{array}\right)
\left(\begin{array}{cc} \bar{\psi}_1 & \bar{\psi}_2 \\ -\psi_2 & \psi_1 \end{array}\right),
$$
where $|r|^2 = \sum_{k=1}^3 (x^k)^2$,
and we derive that, by (\ref{normal}),
$$
K_{11} = \frac{1}{|r|^2} (x^1
 (\psi_1\psi_2 -\bar{\psi}_1\bar{\psi}_2) - ix^2 (\psi_1\psi_2 +\bar{\psi}_1\bar{\psi}_2)+ix^3
(|\psi_2|^2 - |\psi_1|^2)) =
$$
$$
=
\frac{i}{|r|^2} e^\alpha \langle r, \n\rangle.
$$
It follows from (\ref{tildepotential}) that $\widetilde{U} = U + W$ is the potential of $\widetilde{S}$.
Theorem is proved.

More tedious computations allow to prove Theorem \ref{th2} by comparing (\ref{kmatrix}) and (\ref{inversionu}).

For completeness, let us derive (\ref{moutard1}) in the framework of  Section \ref{section4}
(in \cite{C} its derivation was skipped).

Let
$$
\D \Psi = \D \Psi_0 =  0, \ \ \ \widetilde{D} \widetilde{\Psi}_0 = 0,
$$
and let us look for a deformation $\widetilde{\Psi}$ of $\Psi$ of the form
$$
\widetilde{\Psi} = \Psi + \widetilde{\Psi}_0 N
$$
which satisfies $\widetilde{\D}\widetilde{\Psi}=0$.
By (\ref{leibniz}), we have
$$
0 = \widetilde{\D}\widetilde{\Psi} = (\D + W)(\Psi + \widetilde{\Psi}_0 N) = \D\Psi + W\Psi +
(\widetilde{\D}\widetilde{\Psi}_0)\cdot N +
$$
$$
\left(\begin{array}{cc} 0 & 1 \\ 0 & 0 \end{array}\right)\widetilde{\Psi}_0 \partial N +
\left(\begin{array}{cc} 0 & 0 \\ -1 & 0 \end{array}\right)\widetilde{\Psi}_0 \bar{\partial} N
$$
where $W = \widetilde{U}-U$
and, since  $\widetilde{\D}\widetilde{\Psi} = \D \Psi=0$, we look for $N$ such that
$$
W\Psi = -\left(\begin{array}{cc} 0 & 1 \\ 0 & 0 \end{array}\right)\widetilde{\Psi}_0 \partial N -
\left(\begin{array}{cc} 0 & 0 \\ -1 & 0 \end{array}\right)\widetilde{\Psi}_0 \bar{\partial} N.
$$
Recall that, by (\ref{potential}),
$$
W = \left(\begin{array}{cc} 0 & 1 \\ 0 & 0 \end{array}\right)\widetilde{\Psi}_0 S_z (\Psi_0,\Psi)\Psi^{-1} +
\left(\begin{array}{cc} 0 & 0 \\ -1 & 0 \end{array}\right)\widetilde{\Psi}_0 S_{\bar{z}}(\Psi_o,\Psi) \Psi^{-1}
$$
and infer the following

\begin{proposition}[\cite{C}]
\label{prop5}
If $\Psi$ satisfies $\D\Psi =0$, then the function
\begin{equation}
\label{floquet}
\widetilde{\Psi} = \Psi - \widetilde{\Psi}_0 S(\Psi_0,\Psi)
\end{equation}
satisfies the equation $\widetilde{\D}\widetilde{\Psi} = 0$. The function $\widetilde{\Psi}$ is defined
up to $\widetilde{\Psi}_0 \cdot A$ where $A$ is a constant matrix from $H$.
\end{proposition}

Let the potential $U$ is double-periodic:
$$
U(z+\lambda)=U(z), \ \ \ \lambda \in \Lambda \approx {\mathbb Z}^2 \subset \C.
$$
The a solution $\psi$ of the Dirac equation (\ref{diraceq0}) is called the Floquet function (on the zero energy level)
of $\D$ if there are constants $\mu_1$ and $\mu_2$ (the Floquet multipliers) such that
$$
\psi(z+\lambda_k) = \mu_k \psi(z), \ \ \ k=1,2,
$$
where $\lambda_1$ and $\lambda_2$ generate $\Lambda$. The Floquet functions are parameterized by the spectral curve
of $\D$ \cite{T98} which was first introduced for the two-dimensional Schr\"odinger operator in \cite{DKN}.

We conjectured that the spectral curve is preserved by conformal transformations. Since translations and rotations
do not change the potential, it is enough to establish that for the inversion. In \cite{GS} this conjecture was confirmed
by proving that for tori in $\R^3$ the multipliers are preserved by the inversion. In \cite{GT} that was established for
tori in $\R^4$ however it was shown that in this case the spectral curve may change by stacking and unstacking multiple points.
Both proofs are based on studying infinitesimal conformal transformations.

Proposition \ref{prop5} straightforwardly shows the multipliers are preserved by the inversion:

\begin{corollary}
If $\Psi$ is the Floquet function of $\D$, then there is a unique choice of $\widetilde{\Psi}$ of the form (\ref{floquet}) such that
$\widetilde{\Psi}$ is a Floquet function of $\widetilde{\D}$. Moreover $\widetilde{\Psi}$ has the same Floquet multipliers $\mu_1, \mu_2$ as $\Psi$.
\end{corollary}

\section{Final remarks}

1) The (Bianchi) permutability theorem for the Moutard transformation is briefly mentioned in \cite{C}.
It reads that there are representatives of $\widetilde{\Psi}_1$ and
$\widetilde{\Psi}_2$ such that $\widehat{U} = U_{12}=U_{21}$, i.e. the diagram
$$
\begin{array}{ccccc} & U & \stackrel{\Psi_1}{\longrightarrow} & U_1 & \\
 \Psi_2  & \downarrow & & \downarrow &  \widetilde{\Psi}_2  \\
& U_2 & \stackrel{\widetilde{\Psi}_1}{\longrightarrow} & \widehat{U} &
\end{array}
$$
is commutative.
Here the arrows denote transformations defined by its labels and
$\widetilde{\Psi}_2$ and $\widetilde{\Psi}_2$ are images of $\Psi_2$ and $\Psi_1$ under these
transformations. Is there any geometrical interpretation of the permutability in the spirit of Theorem \ref{th2}?

2) In \cite{TT} iterates of the Moutard transformation of the Schr\"o\-din\-ger operator were used for deriving examples of operators with interesting spectral properties. It would be interesting to do the same for the Dirac operator, in particular, by using the geometrical interpretation given by Theorem \ref{th2}.

3)  Proposition \ref{prop3} gives a way for looking for blowups of solutions to the mNV equation. Indeed, let the initial surface $S$ do not pass through $\x=0$. Then (\ref{omega1}) and (\ref{surface}) define its deformation and as soon as $\widetilde{S}$ passes through
$\x=0$,  the potential $\widetilde{U}$  would become singular.  The well-understood minimal surfaces and soliton spheres \cite{T-WN,BP} may supply such explicit examples.

4) There have to be a similar Moutard transformation corresponding to the inversion of surfaces in $\R^4$ \cite{K2,TDS}.
We remark that in this case the Dirac operator takes the form $\D_0 + \mathrm{diag}(U,\bar{U})$ and
is related to the Davey--Stewartson equations.

\end{document}